\documentclass[twoside, 11pt]{article}

\usepackage{mathrsfs,amsfonts,amsmath}
\RequirePackage{ifthen,calc}
\usepackage{theorem}
\usepackage[dvips]{color}

\def\R{{\mathbb R}}

\def\E{{\mathbb E}}

\def\N{{\mathbb N}}

 \catcode`@=11
 \def\@evenhead{\hbox to\textwidth{\footnotesize\rm\thepage \hfill
  {\it }}} % authors name

 \def\@oddhead{\hbox to \textwidth{\footnotesize{\it
 } \hfill\thepage}}% abbreviate title

 \renewcommand{\section}{\makeatletter
 \renewcommand{\@seccntformat}[1]{{\csname the##1\endcsname.}\hspace{0.45em}}
 \makeatother \@startsection
{section}%                                            the name
{1}%                                                  the level
{0pt}%                                                the indent
{\baselineskip}%                                      the beforeskip
{0.5\baselineskip}%                                   the afterskip
{\normalsize\bfseries\mathversion{bold}}}

\renewcommand{\subsection}{\makeatletter
 \renewcommand{\@seccntformat}[1]{{\csname the##1\endcsname.}\hspace{0.45em}}
 \makeatother \@startsection
{subsection}%                                            the name
{1}%                                                  the level
{0pt}%                                                the indent
{\baselineskip}%                                      the beforeskip
{0.5\baselineskip}%                                   the afterskip
{\normalsize\bfseries\mathversion{bold}}}
\catcode`@=12

\newtheorem{thm}{\noindent Theorem}[section]
\newtheorem{lem}{\noindent Lemma}[section]

\theoremheaderfont{\normalfont\bfseries}
\theorembodyfont{\slshape} {\theorembodyfont{\rmfamily}} {\theorembodyfont{\rmfamily}
}
{\theorembodyfont{\rmfamily} } {\theorembodyfont{\rmfamily}
}
{\theorembodyfont{\rmfamily} } {\theorembodyfont{\rmfamily}
}
{\theorembodyfont{\rmfamily} } {\theorembodyfont{\rmfamily} } {\theorembodyfont{\rmfamily}
}
 \def\beqlb{\begin{eqnarray}}\def\eeqlb{\end{eqnarray}}
 \def\beqnn{\begin{eqnarray*}}\def\eeqnn{\end{eqnarray*}}
 \numberwithin{equation}{section}
\def\2R{\mathbb{R}_+\times\mathbb{R}}
\def\qed{\hfill$\square$\smallskip}

\begin{document}
\title{\bf   Random Walks and  Subfractional Brownian Motion}
\author{\small Hongshuai Dai \thanks{E-mail:mathdsh@gmail.com} \\\small  School of Mathematics, Guangxi
University, \\\small Nanning, 530004 China.
\\ \small School of Mathematics and Statistics, Carleton University,
\\ \small Ottawa, K1S 5B6 Canada }
\date{}

\maketitle

\begin{abstract}
\noindent

In this  paper we show an approximation  in law to the subfractional Brownian motion with $H>\frac{1}{2}$ in the Skorohod topology.  The
construction of these approximations is based on random walks.



\end{abstract}

\bigskip

\section{Introduction}

The long-range dependence property is  an important aspect of
stochastic models in various scientific areas,
such as hydrology, telecommunication, finance and so on.
The best known and widely used process that has the long-range dependence
property is the fractional Brownian motion (fBm)
introduced by Mandelbrot and Van Ness (1968).
The fBm is a suitable generalization of the standard Brownian motion and has stationary increments. In many applications, fBm seems to fit very well to random phenomena.
Refer to Samorodnitsky  and  Taqqu (1994) for more information on fBm.

Many scholars have proposed to use other self-similar
Gaussian processes or random fields as stochastic models. This induced some recent progress such as  the generalization of fBm. Many new generalized processes of fBms have been obtained in recent years. As  an  extension  of fBm, subfractional Brownian motion (sub-fBm) has been proposed independently by Bojdecki, Gorostiza and Talarczyk (2004) and by Dzhaparidze and Van Zanten (2004). This process arises from occupation time fluctuations of branching particle systems with Poisson initial condition.
Recall that  the subfractional Brownian motion $X^{H}=\{X^{H}(t),t\geq 0\}$ with index $H\in(0,\;1)$ is a centered Gaussian process with $X^H(0)=0$ and the covariance function
\beqnn
\E\big[X^H(t)X^H(s)\big]=s^{2H}+t^{2H}-\frac{1}{2}\big[(s+t)^{2H}+|t-s|^{2H}\big].
\eeqnn
Sub-fBm $X^H$ is neither a semi-martingale nor a Markov process
unless $H=\frac{1}{2}$. When $H=\frac{1}{2}$,
sub-fBm and the standard Brownian motion coincide.
Sub-fBm has properties analogous to those of fBm (self-similarity, long-range dependence and H$\ddot{o}$lder paths) and satisfies the following inequalities for $s<t$
\beqlb\label{2-2}
\big[(2-2^{2H-1})\wedge 1\big](t-s)^{2H}\leq \E\Big[\big(X^H(t)-X^H(s)\big)^2\Big]\leq \big[(2-2^{2H-1})\vee 1\big]^{2H}(t-s)^{2H}.\qquad
\eeqlb However its increments are not stationary. More information on sub-fBm can be found in Tudor (2007, 2009),  Yan and his coauthors (2010, 2011).
Sub-fBm has raised many interesting theoretical questions. However, in contrast to the extensive study on fBm,
there has been little systematic investigation on sub-fBm. The main reason for this is the complexity of dependence structures.
Therefore, it seems interesting to study this process.

Weak convergence to  fBm processes has been studied extensively
since the works of Davydov (1970) and Taqqu (1975).
In recent years,  many new results on approximations of fBms have been established.
See, for example,  Meyer, Sellan and Taqqu (1999), Li and Dai (2011) and the references therein. On the other hand,  weak limit theorems for sub-fBms have  attracted a lot of interest as well.
For example, Bardina and Bascompte (2010) presented a weak theorem for sub-fBm based on a Poisson process.
Harnett and Nualart (2012) proved  weak convergence of some functionals of  sub-fBm.
Garz\'on,  Gorostiza and Le\'on (2012) proved a strong uniform approximation with a rate of convergence
for sub-fBm by means of transport processes. Similar to these works, in this short note we present an approximation to  sub-fBm with $H>\frac{1}{2}$, however based on  random walks.

The rest of  this paper is organized as follows. In Section 2, we recall some preliminaries  and present the main result of this paper.
Section 3 provides the proof of the main result.

\section{Preliminaries and Main Result}

Let $X^H=\{X^H(t), t\geq 0\}$ be the subfractional Brownian motion with index $H\in(\frac{1}{2},\;1)$.
By using the Hankel transform in  Dzhaparidze and Van Zanten (2004),
we get  from  Tudor (2009) that the process
\beqnn
W(t)=\int_0^{t}\phi_{H}(t,s)dX^{H}(s)
\eeqnn
is the unique Brownian motion such that
\beqlb\label{2-4}
X^H(t)=\int_0^t K_H(t,s) W(ds), \; t\in[0,\;1],
\eeqlb
where
\beqnn
\phi_H(t,s)=\frac{s^{H-\frac{1}{2}}}{\Gamma(\frac{3}{2}-H)}\Big[t^{H-\frac{3}{2}}\big(t^2-s^2\big)^{\frac{1}{2}-H}-\big(H-\frac{3}{2}\big)
\int_s^t\big(x^2-s^2\big)^{\frac{1}{2}-H}x^{H-\frac{3}{2}}dx\Big]1_{(0,\;t)}(s),
\eeqnn
and
\beqnn
K_H(t,s)=\frac{C_H\sqrt{\pi}}{2^{H-1}\Gamma(H-\frac{1}{2})}s^{\frac{3}{2}-H}\int_s^t(x^2-s^2)^{H-\frac{3}{2}}dx1_{(0,\;t)}(s)
\eeqnn
with $C_H$  a normalizing constant. Moreover, $X^H$ and $W$ generate the same filtration.

Let us recall some known facts. Consider a sequence $\{\xi_i\}_{i\in\N}$ of I.I.D random variables  with $\E[\xi_i]=0$ and $\E[\xi^2_i]=1$. The Donsker's invariance principle states that the sequence of the processes
\beqnn
W_n(t)=\sum_{k=1}^{\lfloor nt\rfloor}\frac{\xi_k}{\sqrt{n}}
\eeqnn
converges weakly to a Brownian motion in the Skorohod  topology. Here $\lfloor x\rfloor$ stands for the greatest integer not exceeding $x$.

This result has been extended by Sottinen (2001) to  fractional Brownian motion. Define
\beqnn
F_n(t,s)=n\int_{s-\frac{1}{n}}^s F(\frac{\lfloor nt\rfloor}{n},u)du,\;n\geq 1
\eeqnn
where $F$ is the kernel that transforms the standard Brownian motion into a fractional Brownian one, i.e.,
\beqnn
F(t,s)=c_H (H-\frac{1}{2})s^{\frac{1}{2}-H}\int_{s}^t u^{H-\frac{1}{2}}(u-s)^{H-\frac{3}{2}}du,
\eeqnn
where $c_H$ is still a normalizing constant. Set
\beqnn
Z_n(t)=\int_0^t F_n(t,s)W_n(ds)=\sum_{i=1}^{\lfloor nt\rfloor}n\int_{\frac{i-1}{n}}^\frac{i}{n}F(\frac{\lfloor nt\rfloor}{n},u)du\frac{\xi_i}{\sqrt{n}}.
\eeqnn
Then $\{Z_n(t)\}$ converges weakly to fBm with $H>\frac{1}{2}$. In this paper we mainly extend the above result to the sub-fBm $X^H$ with $H>\frac{1}{2}$.  So we assume  $H>\frac{1}{2}$ in the rest of this paper.

Inspired by Sottinen (2001), we define:
\beqnn
K_n(t,s)=n\int_{s-\frac{1}{n}}^s K_H(\frac{\lfloor nt\rfloor}{n},u)du,\;n\geq 1.
\eeqnn
It is obvious that $K_n(t,\cdot)$¡¡ is an approximation of $K_H(t,\cdot)$ for every $t\in[0,\;1]$. Let
\beqlb\label{2-7}
X_n(t)=n\int_{0}^t K_n(t,u) W_n(du)=\sum_{k=1}^{\lfloor nt \rfloor}n\int_{\frac{k-1}{n}}^{\frac{k}{n}}K_H(\frac{\lfloor nt\rfloor}{n},u)du\frac{\xi_k}{\sqrt{n}}.
\eeqlb
In this paper, we will prove
\begin{thm}\label{thm}
The family of processes $\{X_n(t),\;t\in [0,\;1]\}$  converges weakly in the Skorohod toplogy, as $n$ tends to infinity, to the sub-fractional Brownian motion $X^H$ given by (\ref{2-4}).
\end{thm}

In the rest of this paper, most of  estimates  contain unspecified constants. An unspecified positive and finite constant will be denoted by $C$, which may not be the same in each occurrence. Sometimes we shall
emphasize the dependence of these constants upon parameters.
\section{Proof of Theorem \ref{thm}}

In this section, we will prove Theorem \ref{thm}.  In order to reach our aim, we first verify the convergence of finite-dimensional distributions.

\begin{lem}\label{3-lem}
The family of stochastic processes $\{X_n(t),\;t\in[0,\;1]\}$ given by (\ref{2-7}) converges in the sense of finite-dimensional distributions to the sub-fBm $X^H$ defined by (\ref{2-4}).
\end{lem}

\noindent{\it Proof:} In order to prove Lemma \ref{3-lem}, it suffices to prove that for any $t_1,\cdots,t_p\in[0,\;1]$ and $\eta\in\R$, we have
\beqlb\label{a-1}
\E\Bigg[\exp\Big(i\eta\sum_{i=1}^pX_n(t_i)\Big)\Bigg]\to \E\Bigg[\exp\Big(i\eta\sum_{i=1}^pX^H(t_i)\Big)\Bigg]
\eeqlb
as $n\to\infty$.

In order to prove (\ref{a-1}), we first need to introduce the following notation. Let us consider a  sequence  of partitions $\{t^m_i\}$ of the interval $[0,\,1]$ of the form
\beqnn
\pi^m:\;0=t^m_0<t^m_1<t^m_2<\cdots<t^m_{m}=1,
\eeqnn
where $t_i^m=\frac{i}{m}$, $i=0,1,\cdots,m$. We  define

\beqnn
\hat{X}_{m}(t)=\sum_{i=1}^{m}K_H(t,t^m_{i-1})W(\Delta_i),
\eeqnn
and
\beqlb\label{3-3}
X_{m,n}(t)=\sum_{i=1}^{m}K_H(t,t^m_{i-1})W_n(\Delta_i),
\eeqlb
where $\Delta_i=[t^m_{i-1}, \;t^m_i)$,
$
W(\Delta_i)=W(t^m_i)-W(t^m_{i-1}),
$
and $
W_n(\Delta_i)=W_n(t^m_{i})-W_n(t^m_{i-1}).
$
Moreover, let
\beqnn
\tilde{X}_{n}(t)&&=\sum_{i=1}^{\lfloor nt\rfloor}n \int_{\frac{i-1}{n}}^{\frac{i}{n}} K_H(t,u)du \frac{\xi_i}{\sqrt{n}},
\\ \tilde{X}_{m,n}(t)&&=\sum_{k=1}^{\lfloor nt\rfloor}n\int_{\frac{k-1}{n}}^{\frac{k}{n}}\sum_{i=1}^{m}K_H(t,t^m_{i-1})1_{\Delta_i}(u)du\frac{\xi_k}{\sqrt{n}},
\eeqnn
and
\beqnn
K^m(t,s)=\sum_{i=1}^{m}K_H(t,t^m_{i-1})1_{\Delta_i}(s).
\eeqnn
One can easily get that
\beqnn
K^m(t,\cdot)\to K_H(t,\cdot)
\eeqnn
as $m\to\infty$ in $L^2\big([0,\;1]\big)$.

We can easily get that
\beqlb\label{a-22}
&&\Bigg|\E\bigg[\exp\Big(i\eta\sum_{i=1}^pX_n(t_i)\Big)-\exp\Big(i\eta\sum_{i=1}^pX^H(t_i)\Big)\bigg]\Bigg|\nonumber
\\&&\qquad\qquad\leq |D_1(n)|+|D_2(n,m)|+|D_3(n,m)|+|D_4(n,m)|+|D_5(m)|,\qquad
\eeqlb
where
\beqnn
D_1(n)&&=\E\bigg[\exp\Big(i\eta\sum_{i=1}^pX_n(t_i)\Big)-\exp\Big(i\eta\sum_{i=1}^p\tilde{X}_n(t_i)\Big)\bigg],
\\
D_2(n,m)&&=\E\bigg[\exp\Big(i\eta\sum_{i=1}^p\tilde{X}_{m,n}(t_i)\Big)-\exp\Big(i\eta\sum_{i=1}^p\tilde{X}_n(t_i)\Big)\bigg],
\\
D_3(n,m)&&=\E\bigg[\exp\Big(i\eta\sum_{i=1}^p\tilde{X}_{m,n}(t_i)\Big)-\exp\Big(i\eta\sum_{i=1}^pX_{m,n}(t_i)\Big)\bigg],
\\
D_4(n,m)&&=\E\bigg[\exp\Big(i\eta\sum_{i=1}^pX_{m,n}(t_i)\Big)-\exp\Big(i\eta\sum_{i=1}^p\hat{X}_m(t_i)\Big)\bigg],
\eeqnn
and
\beqnn
D_5(m)&&=\E\bigg[\exp\Big(i\eta\sum_{i=1}^p \hat{X}_{m}(t_i)\Big)-\exp\Big(i\eta\sum_{i=1}^pX^H(t_i)\Big)\bigg].
\eeqnn
 In order to simplify our discussion, here we assume that $p=1$. For  $p> 1$, we can use the same method to get the result.  For convenience,  let $t_1=t$.

We first study $D_1(n)$. Since $\xi_i$, $i=1,\cdots,$ are I.I.D, we have that as $n\to\infty$,

\beqnn
\E\Big[\tilde{X}_n(t)-X_n(t)\Big]^2&&=\E\Big[\sum_{k=1}^{\lfloor nt \rfloor}n\int_{\frac{k-1}{n}}^{\frac{k}{n}}\big(K_H(t,u)-K_H(\frac{\lfloor nt \rfloor}{n},u)\big)du\frac{\xi_k}{\sqrt{n}}\Big]^2\nonumber
\\&&=\sum_{k=1}^{\lfloor nt \rfloor}n\Bigg[\int_{\frac{k-1}{n}}^{\frac{k}{n}}\big(K_H(t,u)-K_H(\frac{\lfloor nt \rfloor}{n},u)\big)du\Bigg]^2\nonumber
\\&&\leq C \int_{0}^1\big(K_H(t,u)-K_H(\frac{\lfloor nt \rfloor}{n},u)\big)^2du
\\&&\leq C \Big(t-\frac{\lfloor nt \rfloor}{n}\Big)^{2H}\to 0,
\eeqnn
where we used (\ref{2-2}).
Therefore, as $n\to\infty$
\beqlb\label{3-18}
\E\Big[X_n(t)-\tilde{X}_n(t)\Big]^2\to 0.
\eeqlb
From (\ref{3-18}), one can easily get that as $n\to\infty$,
\beqlb\label{3-23}
|D_1(n)|\to 0.
\eeqlb

Next, we deal with $D_2(n,m)$.  Indeed,  by the fact that
\beqnn
\E[\xi_i\xi_j]=\left\{\begin{array}{lll}0,\;&\textrm{if}\; i\neq j,
\\1,&\textrm{if}\;i=j,\end{array}\right.\eeqnn and the Cauchy-Schwartz inequality, we obtain
\beqnn
\E\Big[\tilde{X}_{m,n}(t)-\tilde{X}_n(t)\Big]^2&&=\E\Big[\sum_{k=1}^{\lfloor nt\rfloor}\Big(n\int_{\frac{k-1}{n}}^{\frac{k}{n}}\big(K^m(t,u)-K_H(t,u)\big)du\Big)\frac{\xi_k}{\sqrt{n}}\Big]^2
\\&&=\sum_{k=1}^{\lfloor nt\rfloor}n\Big(\int_{\frac{k-1}{n}}^{\frac{k}{n}}\big(K^m(t,u)-K_H(t,u)\big)du\Big)^2
\\&&\leq \sum_{k=1}^{\lfloor nt\rfloor}\int_{\frac{k-1}{n}}^{\frac{k}{n}}\big(K^m(t,u)-K_H(t,u)\big)^2du
\\&&\leq \int_0^1\big(K^m(t,u)-K_H(t,u)\big)^2du \to 0
\eeqnn
as $m\to\infty$. Therefore,
 \beqnn
\tilde{X}_{m,n}(t)\to \tilde{X}_{n}(t)
\eeqnn
as $m\to\infty$ with respect to $n$ uniformly in $L^2(\Omega)$.    So
\beqlb\label{3-24}
|D_2(n,m)|\to 0
\eeqlb
as $m\to\infty$ uniformly in $n$.

Now, we deal with $D_3(n,m)$.  Let $\lambda$ be the Lebesque measure.  We also note that if $t^m_i>t$, then  $K_H(t,t_{i}^m)=0$.  Noting that $\int_{\frac{i-1}{n}}^{\frac{i}{n}}1_{\Delta_k}(u)du=\lambda\big([\frac{i-1}{n},\frac{i}{n})\cap
\Delta_k\big)$,
 we can rewrite $\tilde{X}_{m,n}$  as follows,

\beqlb\label{3-9}
\tilde{X}_{m,n}(t)&&=\sum_{i=1}^{m}K_H(t,t^m_{i-1})\sum_{k=1}^{\lfloor nt\rfloor}n\lambda\Big(\big[\frac{k-1}{n}, \frac{k}{n}\big)\cap\Delta_i\Big)\frac{\xi_k}{\sqrt{n}}.
\eeqlb
On the other hand, the points $\frac{i}{n}$, $i=0,\cdots,n,$ also form a  partition of the interval $[0,\;1]$. We let $n$ be sufficiently large and then note that
$$
[\frac{\lfloor nt_{i-1}^m\rfloor+1}{n},\;\frac{\lfloor nt_{i}^m\rfloor}{n})\subset\Delta_i=[t_{i-1}^{m},\;t_i^m).
$$
By discussing the relation between the endpoints $t_{i-1}^m$, $t_i^m$ and the intervals $[\frac{\lfloor nt_{i-1}^m\rfloor}{n},\;\frac{\lfloor nt_{i-1}^m\rfloor+1}{n})$ and $[\frac{\lfloor nt_{i}^m\rfloor}{n},\;\frac{\lfloor nt_{i}^m\rfloor+1}{n})$, we can rewrite \eqref{3-9} as
\beqlb\label{3-9-2}
\tilde{X}_{m,n}(t)=\sum_{i=1}^{m}K_H(t,t^m_{i-1})&&\bigg(\sum_{k=\lfloor nt_{i-1}^m\rfloor+1}^{\lfloor nt_i^m\rfloor}\frac{\xi_k}{\sqrt{n}}-\big(nt_{i-1}^m-\lfloor nt_{i-1}^m\rfloor\big)\frac{\xi_{\lfloor nt_{i-1}^m\rfloor+1}}{\sqrt{n}}\nonumber
\\&& \qquad +\big(nt_i^m-\lfloor nt_i^m\rfloor\big)\frac{\xi_{\lfloor nt_{i}^m\rfloor+1}}{\sqrt{n}}\bigg).
\eeqlb
From (\ref{3-3}), we can get
\beqlb\label{3-8}
X_{m,n}(t)&&=\sum_{i=1}^{m}K_H(t,t^m_{i-1})\sum_{k=\lfloor nt^m_{i-1}\rfloor+1}^{\lfloor nt^m_{i}\rfloor}\frac{\xi_k}{\sqrt{n}}.
\eeqlb
Since $0\leq nt_{i-1}^m-\lfloor nt_{i-1}^m\rfloor \leq 1$ and $0\leq nt_i^m-\lfloor nt_i^m\rfloor\leq 1$, by (\ref{3-9-2}) and (\ref{3-8}),  we can get that
\beqnn
\E\big(|X_{m,n}(t)-\tilde{X}_{m,n}(t)|^2\big)\leq \frac{C}{n}\sum_{i=1}^m K_H^2(t,t_{i-1}^m),
\eeqnn
which implies that for any given $m$,
\beqlb\label{3-25}
|D_3(n,m)|\to 0
\eeqlb
as $n\to \infty$.

Below, we deal with $D_4(n,m)$.   By the invariance principle and the continuous mapping theorem (see e.g.,  Billingsley,  1968), one can easily get that for any given $m$,
\beqlb\label{3-4}
X_{m,n}(t) \stackrel{W}{\Rightarrow}\hat{X}_{m}(t)
\eeqlb
as $n\to \infty$, where $\stackrel{W}{\Rightarrow}$ denotes weak convergence.
By (\ref{3-4}), one can easily get that for any given $m$,
\beqlb\label{3-27}
|D_4(n,m)|\to 0
\eeqlb
as $n\to\infty$.

 Finally, we study $D_5(m)$.  Observing that $K_H(t,s)$ is continuous  in $s$ for every $t$, we can easily get  that
\beqnn
\hat{X}_m(t)\to X^H(t)
\eeqnn
 as $m\to\infty$  in $L^2(\Omega)$. Therefore,
 \beqlb\label{3-26}
 |D_5(m)|\to 0
 \eeqlb
as $m\to\infty$.

By (\ref{a-22}), (\ref{3-23}), (\ref{3-24}), (\ref{3-25}), (\ref{3-27}) and (\ref{3-26}), we can get that (\ref{a-1}) holds.
\qed

Next, we prove the tightness of $\{X_n(t)\}_{n\in\N}$.

\begin{lem}\label{3-lem2}
The family $\{X_n(t),t\in[0,\;1]\}$ given by (\ref{2-7}) is tight.
\end{lem}

\noindent{\it Proof:} Noting that  the kernel $K_H(t,u)$  vanishes when $u$ is larger than $t$, we have that for any $t>s$,
\beqlb\label{3-19}
\E\bigg[X_n(t)-X_n(s)\bigg]^2&&=\E\bigg[\sum_{i=1}^{\lfloor nt\rfloor}n\int_{\frac{i-1}{n}}^{\frac{i}{n}}\Big(K_H(\frac{\lfloor nt\rfloor}{n},u)-K_H(\frac{\lfloor ns\rfloor}{n},u)\Big)du\frac{\xi_i}{\sqrt{n}}\bigg]^2\nonumber
\\&&=\sum_{i=1}^{\lfloor nt\rfloor}n\Big(\int_{\frac{i-1}{n}}^{\frac{i}{n}}\Big(K_H(\frac{\lfloor nt\rfloor}{n},u)-K_H(\frac{\lfloor ns\rfloor}{n},u)\Big)du\Big)^2,
\eeqlb
since $\E[\xi_i\xi_j]=0$ if $i\neq j$, and $\E[\xi_i^2]=1$.

By the H$\ddot{o}$lder inequality,  (\ref{3-19}) can be bounded by
\beqlb\label{3-20}
\sum_{i=1}^{\lfloor nt\rfloor}\int_{\frac{i-1}{n}}^{\frac{i}{n}}\Big(K_H(\frac{\lfloor nt\rfloor}{n},u)-K_H(\frac{\lfloor ns\rfloor}{n},u)\Big)^2du&&\leq C\int_0^1\Big(K_H(\frac{\lfloor nt\rfloor}{n},u)-K_H(\frac{\lfloor ns\rfloor}{n},u)\Big)^2du\nonumber
\\&&\leq C\Big|\frac{\lfloor nt\rfloor}{n}-\frac{\lfloor ns\rfloor}{n}\Big|^{2H},
\eeqlb
since (\ref{2-2}) holds. By using the same argument as in Torres and Tudor (2009), we can get from \eqref{3-20} that for any $s<t<u\in[0,\;1]$,
\beqlb\label{3-a1}
\E|X_n(t)-X_n(s)||X_n(u)-X_n(t)|\leq C|u-s|^{2H}.
\eeqlb
By (\ref{3-a1}) and Billingsley  (1968), we can get that the lemma holds, since $H>\frac{1}{2}$.
\qed

Now, we prove the main result of this paper.

\noindent{\it Proof of Theorem \ref{thm}:}  Theorem \ref{thm} is a direct consequence of Lemmas \ref{3-lem} and \ref{3-lem2}, because tightness and the convergence of finite dimensional distributions imply  weak convergence in the Skorohod topology (see Billingsly, 1968). \qed

\noindent{\bf Acknowledgments}\ The author thanks two referees for very detailed comments and suggestions which resulted in  short proofs and better presentation of this paper. This work was  supported by the National Natural Science Foundation of China (11361007, 11061002) and the Guangxi Natural Science Foundation (2012GXNSFBA053010, 2011GXNSFA018126).


\begin{thebibliography}{hhhh}
\bibliographystyle{100}
\bibitem{B1968}
Billingsley, P.~(1968).
\newblock{\em Convergence of Probability Measures.}
\newblock{ New York: John Willey and Sons.}
\bibitem{BB}
Bardina, B., Bascompte, D.~(2010).
\newblock{Weak convergence towards  two independent Gaussian processes from a unique poisson process.}
\newblock{\em Collect. Math.} 61:191-204.
\bibitem{BGT}
Bojdecki, T., Gorostiza, L.~G., Talarczyk, A.~(2004).
\newblock{Sub-fractional Brownian motion and its relation to occupation times. }
\newblock{\em Statist. Probab. Lett.} 69: 405-419.


\bibitem{D70}
Davydov, Y.~(1970).
\newblock{The invariance principle for stationary processes.}
\newblock{\em Teor. Veroatn. Ee Primenen.} 15: 498-509.
\bibitem{DVZ}
Dzhaparidze, K., Van Zanten, H.~(2004).
\newblock{A series expansion of fractional Brownian motion.}
\newblock{\em Probab. Theory Relat. Fields} 103: 39-55.


\bibitem{GG}
Garz\'on, J., Gorostiza, G., Le\'on, A.~(2012).
\newblock{A strong uniform approximation of sub-fractional Brownian motion.}
\newblock{Preprint.}
\bibitem{HN}
Harnett, D., Nualart, D.~(2012).
\newblock{Weak convergence of the stratonovich integral with respect to a class of Gaussian processes.}
\newblock{\em Stochastic Processes and Their Applications} 122: 3460-3505.
\bibitem{LD2011}
Li, Y., Dai, H.~(2011).
\newblock{Approximations of fractional Brownian motion.}
\newblock{\em Bernoulli} 17: 1195-1216.
\bibitem{MV1968}
Mandelbrot, B., Van Ness, J.~W.~(1968).
\newblock{Fractional Brownian motions, fractional noises and applicaitons.}
\newblock{\em SIAM Rev.} 10: 422-437.
\bibitem{MST}
Meyer, Y., Sellan, F., Taqqu, M.~S.~(1999).
\newblock{Wavelets, generalized
white noise and fractional integration: The synthesis of fractional
Brownian motion.}
\newblock {\em J. Fourier Anal. Appl.} 5: 465-494.
\bibitem{ST94}
Samorodnitsky, G., Taqqu, M.~S.~(1994).
\newblock{\em Stable Non-Gaussian Random Processes.}
 New York: Chapman and Hall.
\bibitem{T2001}
Sottinen, T.~(2001).
\newblock{Fractional Brownian motion, random walks and binary market models.}
\newblock{\em Finance and Stochastic} 5: 343-355.
\bibitem{Taqqu75}
Taqqu, M.~S.~(1975).
\newblock{Weak convergence to fractional Brownian
motion and to the Rosenblatt process.}
\newblock{\em Z. Wahrsch. Verw.
Gebiete} 31: 287-302.
\bibitem{TT2009}
Torres, S.,  Tudor, C.~A.~(2009).
\newblock{Donsker type theorem for the Rosenblatt process and a binary market model.}
\newblock{\em Stochastic Analysis and Applications} 27: 555-573.
\bibitem{T2007}
Tudor, C.~A.~(2007).
\newblock{Some properties of the sub-fractional Brownian motion.}
\newblock{\em Stochastics} 79: 431-448.

\bibitem{T2009}
 Tudor, C.~A.~(2009).
\newblock{On the Wiener integral with respect to a sub-fractional Brownian motion on an interval.}
\newblock{\em J. Math. Anal. Appl.} 351: 456-468.

\bibitem{YS2010}
Yan, L., Shen, G.~(2010).
\newblock{On the collision local time of sub-fractional Brownian motions.}
\newblock{\em Statist. Probab. Lett.} 80: 296-308.
\bibitem{YS2011}
Yan, L., Shen, G., He, K.~(2011).
\newblock{It$\hat{o}$'s formula for a subfractional Brownian motion.}
\newblock{\em Commun. Stoch. Anal.} 5: 135-159.
\end{thebibliography}
\end{document}